\documentclass[10pt]{article}

\makeatletter

\newif\iftwodates
\twodatesfalse

\renewcommand\maketitle{\begin{titlepage}%
  \pagenumbering{Alph}  
  \let\footnotesize\small
  \let\footnoterule\relax
  \let\footnote\thanks
  \null\vfil
  \vskip 30\p@
  \begin{center}%
    {\LARGE \bf \@title \par}%
    \vskip 3em%
    {\large
     \lineskip .75em%
     \begin{tabular}[t]{c}%
       \@author
     \end{tabular}\par}%
     \vskip 1.5em%
  \end{center}\par
  \vfill
  \begin{center}
    \raisebox{1.5cm}{\includegraphics[width=0.58\textwidth]%
      {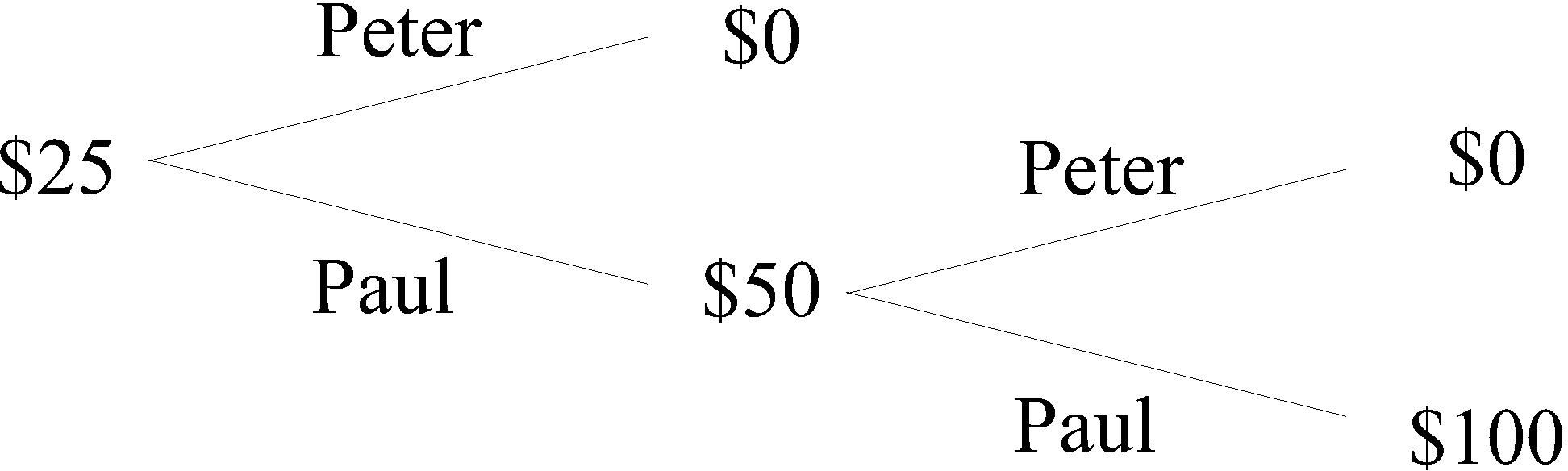}}%
    \hskip 3em%
    \includegraphics[width=0.29\textwidth]%
      {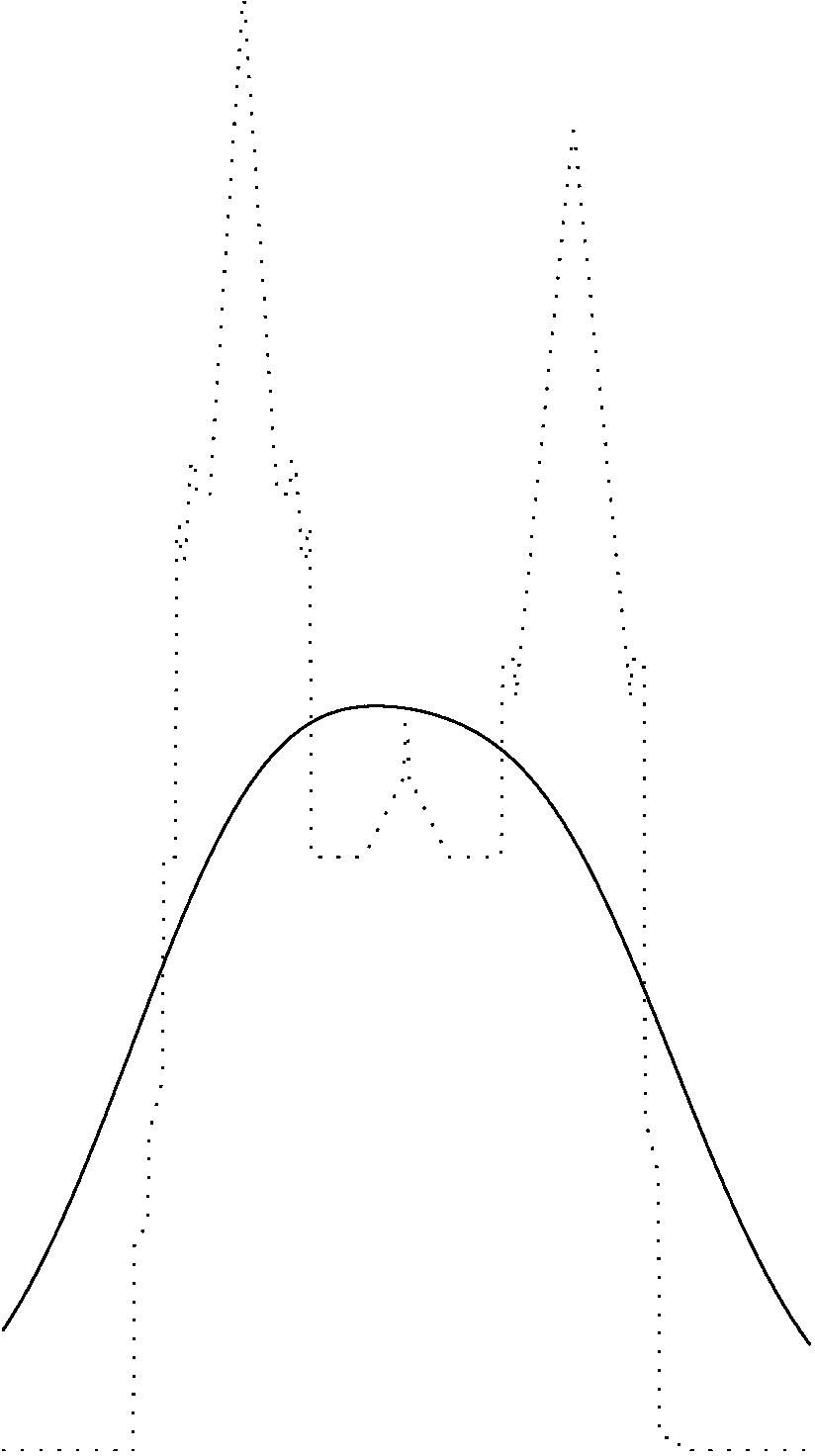}%
  \end{center}
  \@thanks
  \vfill
  \begin{center}
    {\large \bf The Game-Theoretic Probability and Finance Project}
  \end{center}
  \begin{center}
    {\large Working Paper \#\No}
  \end{center}
  \begin{center}
    {\iftwodates\large First posted \firstposted.
    Last revised \@date.\else\large\@date\fi}
  \end{center}
  \begin{center}
    Project web site:\\
    http://www.probabilityandfinance.com
  \end{center}
  \end{titlepage}%
  \setcounter{footnote}{0}%
  \global\let\thanks\relax
  \global\let\maketitle\relax
  \global\let\@thanks\@empty
  \global\let\@author\@empty
  \global\let\@date\@empty
  \global\let\@title\@empty
  \global\let\title\relax
  \global\let\author\relax
  \global\let\date\relax
  \global\let\and\relax
}

\renewenvironment{abstract}{%
  \titlepage\pagenumbering{roman}  
  \null\vfil
  \@beginparpenalty\@lowpenalty
  \begin{center}%
    \Large \bfseries \abstractname
    \@endparpenalty\@M
  \end{center}}%
  {\par\vfill\tableofcontents\thispagestyle{empty}\endtitlepage
  \pagenumbering{arabic}}  

\makeatother
\newcommand*{\No}{54}
\twodatestrue
\newcommand*{\firstposted}{March 5, 2019}

\usepackage{amsmath,amsthm,amssymb,amsfonts,latexsym,graphicx,enumitem,booktabs}
\usepackage[numbers]{natbib}

\usepackage{mathtools}

\usepackage{mathcomp}

\usepackage{array}
\newcommand{\PreserveBackslash}[1]{\let\temp=\\#1\let\\=\temp}
\newcolumntype{C}[1]{>{\PreserveBackslash\centering}p{#1}}
\newcolumntype{R}[1]{>{\PreserveBackslash\raggedleft}p{#1}}
\newcolumntype{L}[1]{>{\PreserveBackslash\raggedright}p{#1}}

\usepackage[utf8]{inputenc}

\makeatletter
\renewcommand\tableofcontents{%
    \@starttoc{toc}%
}
\makeatother

\allowdisplaybreaks[4]

\newif\ifFULL
\FULLfalse

\ifFULL
  \usepackage{color}

\fi

\theoremstyle{definition}

\newlength{\IndentI}
\newlength{\IndentII}
\newlength{\IndentIII}
\newlength{\IndentIV}
\newcommand{\indentI}{\noindent\hspace*{\IndentI}}
\newcommand{\indentII}{\noindent\hspace*{\IndentII}}

\newcommand{\bbbr}{\mathbb{R}}  
\newcommand{\K}{\mathcal{K}}    

\newcommand{\st}{\mathbin{|}}

\title{The language of betting as a strategy for statistical and scientific communication}
\date{\today}

\author{Glenn Shafer, Rutgers University}

\begin{document}


\begin{abstract}
The most widely used concept of statistical inference --- the p-value --- is too complicated for effective communication to a wide audience.  This paper proposes a simplification:  replace p-values with outcomes of bets.  This leads to a new understanding of likelihood, to alternatives to power and confidence, and to a framework that accommodates both planned and opportunistic testing of statistical hypotheses and other probabilistic forecasting systems.  This framework for statistics builds on the framework for mathematical probability developed in \emph{Game-Theoretic Foundations for Probability and Finance} (Glenn Shafer and Vladimir Vovk, Wiley, 2019).
\end{abstract}

\section{Introduction}

The most widely used concept of statistical inference --- the p-value --- is too complicated for effective communication to a wide audience \cite{McShane/Gal:2017,Gigerenzer:2018}.  This paper proposes a simplification:  replace p-values with outcomes of bets.  This leads to a new understanding of likelihood, to alternatives to power and confidence, and to a framework that accommodates both planned and opportunistic testing of statistical hypotheses and other probabilistic forecasting systems. 

Testing a hypothesized probability distribution by betting is straightforward.  We bet by selecting a payoff and buying it for its hypothesized expected value.  If the bet multiplies the money it risks by a large factor, we have evidence against the hypothesis, and the factor measures the strength of this evidence.  Multiplying our money by 5 might merit attention; multiplying it by 100 or by 1000 might be considered conclusive.  

The factor by which we multiply the money we risk --- we may call it the \emph{betting score} --- is conceptually simpler than a p-value, because it reports the result of a single bet, whereas a p-value is based on a family of tests. As explained in Section~\ref{sec:test}, betting scores also have a number of other advantages:
\begin{itemize}
\item
Whereas the certainty provided by a p-value is sometimes exaggerated, the uncertainty remaining when a large betting score is obtained cannot be hidden.  Everyone knows that a bet can succeed by sheer luck.
\item
A bet (a payoff selected and bought) determines an implied alternative hypothesis, and the betting score is the likelihood ratio with respect to this alternative.  So the evidential meaning of betting scores is aligned with our intuitions about likelihood ratios. 
\item
Along with its implied alternative hypothesis, a bet determines an implied target:  a value for the betting score that can be hoped for under the alternative hypothesis.   Implied targets can be more useful than power calculations, because an implied target along with an actual betting score tells a coherent story.  The notion of power, because it requires a fixed significance level, does not similarly cohere with the notion of a p-value.
\item
Testing by betting permits opportunistic searches for significance, because the persuasiveness of having multiplied money risked does not depend on having followed a complete betting strategy laid out in advance.  
\end{itemize}

The message of this paper is that p-values can and should be replaced by betting scores in most statistical practice.  But this cannot happen overnight, and the notion of calculating a p-value will always be on the table when statisticians look at the standard or probable error of the estimate of a difference.%
\footnote{See \cite{Shafer:2019} for a discussion of how this happened in the 1830s.}
We will want, therefore, to relate the scale for measuring evidence provided by a p-value to the scale provided by a betting score.  Any rule for translating from the one scale to the other will be arbitrary, but it may nevertheless be useful to establish some such rule as a standard.  This issue is discussed in Section~\ref{sec:comp}.

Section~\ref{sec:estimation} considers statistical modeling and estimation.  Here betting moves partly off-stage.  Just as the usual picture of a statistical model involves partial knowledge of a probability distribution, the betting picture involves seeing only some moves of a betting game.  We do not see how a bet pays off.  But we can nevertheless equate the model's validity with the futility of betting against it, and a strategy for a hypothetical bettor inside the game, together with observations of the phenomenon, then translates into a \emph{warranty} about one or more unknown parameters.  Instead of ($1-\alpha$)-confidence, we obtain a $(1/\alpha)$-warranty, which tells us that the hypothetical bettor has multiplied the money he risked by $1/\alpha$ if the parameter is not in a particular subset of the parameter space.

A ($1-\alpha$)-confidence set qualifies as a $(1/\alpha)$-warranty set; it arises when all-or-nothing bets are used to decide whether to include parameter values in the warranty set.  But the more general concept of warranty obtained by allowing bets that are not all-or-nothing has several advantages:
\begin{itemize}
\item
It shares with individual betting scores the advantage that its residual uncertainty is not hidden.  
\item
The betting strategy and observations produce more than one warranty set.  They produce a $(1/\alpha)$-warranty set for every $\alpha$, and these warranty sets are nested.
\item
Because it is always legitimate to continue betting with whatever capital remains, the hypothetical bettor can continue betting on additional outcomes, and we can update our warranty sets accordingly without being accused of ``sampling to a foregone conclusion''.  The same principles authorize us to combine warranty sets based on successive studies.
\end{itemize} 

The conclusion of the paper (Section~\ref{sec:conclusion}) summarizes the advantages of testing by betting.  An appendix (Section~\ref{sec:sit}) situates the idea in the broader landscape of theoretical statistics and other proposed remedies for the misunderstanding and misuse of p-values and significance testing.  

For further discussion of betting as a foundation for mathematical probability, statistics, and finance, see \cite{Shafer/Vovk:2019} and related working papers at www.probabilityandfinance.com.  This paper draws on some of the mathematical results reported in Chapter~10 of \cite{Shafer/Vovk:2019}, but the crucial concepts of implied alternative, implied target, and warranty are newly introduced here.

\section{Testing by betting}\label{sec:test}

You claim that a probability distribution $P$ describes a certain phenomenon $Y$.  How can you give content to your claim, and how can I challenge it? 

Assuming that we will later see $Y$'s actual value $y$, a natural way to proceed is to interpret $P$ as a collection of betting offers.   You offer to sell me any payoff $S(Y)$ for its expected value, $\mathbf{E}_P(S)$.  I choose a nonnegative payoff $S$, so that $\mathbf{E}_P(S)$ is all  I risk. Let us call $S$ my \emph{bet}, and let us call the factor by which I multiply the money I risk,
\[ 
      \frac{S(y)}{\mathbf{E}_P(S)},
\] 
my \emph{betting score}.  This score does not change when $S$ is multiplied by a positive constant.  I will usually assume, for simplicity, that $\mathbf{E}_P(S)=1$ and hence that the score is simply $S(y)$.

A large betting score is the best evidence I can have against $P$.  I  have bet against $P$ and won.  On the other hand, the possibility that I was merely lucky remains stubbornly in everyone's view.  By using the language of betting, I have accepted the uncertainty involved in my test and made sure that everyone else is aware of it as well.

I need not risk a lot of money.  I can risk as little as I like --- so little that I am indifferent to losing it and to winning any amount the bet might yield.  So this use of the language of betting is not a chapter in decision theory.  It involves neither the evaluation of utilities nor any Bayesian reasoning.  I am betting merely to make a point.  But whether I use real money or play money, I must declare my bet before the outcome $y$ is revealed, in the situation in which you asserted $P$.  

This section explains how testing by betting can bring greater flexibility and clarity into statistical testing.  Section~\ref{subsec:basic} explains how betting can be more opportunistic than conventional significance testing.  Section~\ref{subsec:likely} explains that a bet implies an alternative hypothesis, and that the betting score is the likelihood ratio with respect to this alternative.  Section~\ref{subsec:target} explains how the alternative hypothesis in turn implies a target for the bet.  Finally, Section~\ref{subsec:examples} uses three simple but representative examples to show how the concepts of betting score and implied target provide a clear and consistent message about the result of a test, in contrast to the confusion that can arise when we use the concepts of p-value and power.

\subsection{Basic advantages}\label{subsec:basic}

The standard way of testing a probability distribution $P$ is to select a \emph{significance level} $\alpha\in(0,1)$, usually small, and a set $E$ of possible values of $Y$ such that $P(Y\in E)=\alpha$.  The event $E$ is the \emph{rejection region}. The probability distribution $P$ is discredited (or \emph{rejected}) if the actual value $y$ is in $E$.  

Although textbooks seldom make the idea explicit, a standard test is often thought of as a bet:  I pay $\$1$ for the payoff $\$S_E$ defined by
\begin{equation}\label{eq:NP}
          S_E(y):=
          \begin{cases}
            \frac{1}{\alpha} & \text{if $y\in E$}\\
            0                       & \text{if $y\not\in E$}.
          \end{cases}
\end{equation}
If $E$ happens, I have multiplied the $\$1$ I risked by $1/\alpha$.  This makes standard testing a special case of testing by betting, the special case where the bet is \emph{all-or-nothing}.  In return for $\$1$, I get either $\$(1/\alpha)$ or $\$0$.

Although statisticians are accustomed to all-or-nothing bets, there are two good reasons for generalizing beyond them.  First, the betting score $S(y)$ from a more general bet is a graduated appraisal of the strength of the evidence against $P$.  Second, when we allow more general bets, testing can be opportunistic.

\paragraph{A betting outcome is a graduated appraisal of evidence.}  

A betting score $S(y)$ appraises the evidence against $P$.  The larger $S(y)$, the stronger the evidence.

A p-value also appraises the evidence against $P$; the smaller the p-value, the stronger the evidence.  But p-values are more complicated than betting scores; they involve a large class, ideally a continuum, of all-or-nothing tests.  To obtain a p-value, we usually begin with function $T$ of $Y$, called a \emph{test statistic}.  In the ideal case, there exists for each significance level $\alpha\in(0,1)$ a number $t_{\alpha}$ such that 
\begin{equation}\label{eq:p1}
         P\left(T\ge t_{\alpha}\right) = \alpha.
\end{equation}
So we have an all-or-nothing test for each $\alpha$:  reject $P$ if $T(y)\ge t_{\alpha}$.  The \emph{p-value}, say $\mathsf{p}(y)$, is the smallest $\alpha$ for which the test rejects:
\begin{equation}\label{eq:p2}
      \mathsf{p}(y) := \inf \{\alpha \st T(y)\ge t_{\alpha}\} = P(T\ge T(y)).
\end{equation}
The larger $T(y)$, the smaller $\mathsf{p}(y)$.

Large values of $T(y)$ are supposed to discredit $P$.  The p-value $\mathsf{p}(y)$ locates the degree of discredit on a scale from zero to one.  But what does the scale mean?  For a mathematical statistician, this question is answered by \eqref{eq:p1} and \eqref{eq:p2}.  For less sophisticated users of statistics, the import of these equations can be elusive.  It is especially elusive for anyone who relies on betting intuitions to understand probability.  \emph{Had I known $y$ in advance}, I could have multiplied my money by $1/\mathsf{p}(y)$ by making an all-or-nothing bet with significance level $\mathsf{p}(y)$.  But I did not know $y$ in advance, and pretending that I did would be cheating.

\paragraph{Betting can be opportunistic.}

The probabilistic predictions that can be associated with a scientific hypothesis usually go beyond a single comprehensive probability distribution.  A scientist might  begin with a joint probability distribution $P$ for a sequence of variables $Y_1,\ldots,Y_N$ and formulate a plan for successive experiments that will allow her to observe them.  But the scientific enterprise is usually more opportunistic.  A scientist might perform an experiment that produces $Y_1$'s value $y_1$ and then decide whether it is worthwhile to perform the further experiment that would produce $Y_2$'s value $y_1$.  Perhaps no one even thought about $Y_2$ at the outset.  One scientist or team tested the hypothesis using $Y_1$, and then, perhaps because the result was promising but not conclusive, some other scientist or team came up with the idea of further testing the hypothesis with a second variable $Y_2$ from a hitherto uncontemplated new experiment or database.

Testing by betting can accommodate this opportunistic nature of scientific research.  Imagine, for example, that I doubt the validity of the probability forecasts made by a particular weather forecaster.  Imagine further that the forecaster decides each day, on a whim, what to forecast that day; perhaps he will give a probability distribution for the amount of rain, perhaps a probability distribution for the temperature at 10:00 a.m., etc.   In spite of his unpredictability, I can try to show that he is a poor forecaster by betting against him.  I start with $\$1$, and each day I buy a random variable for the expected value he attributes to it.  I take care never to risk more than I have accumulated so far, so that my overall risk never exceeds the $\$1$ with which I began.  If I have accumulated $\$1000$ after a year or two, this will be convincing evidence against the forecaster.

Such opportunistic betting boils down to multiplying betting scores.  My initial capital is $1$.  My first bet $S_1$ is nonnegative and has expected value $1$ according to the forecaster.  After it is settled, my capital is $S_1(y_1)$.  Now I select a nonnegative bet $S_2$ to which the forecaster now gives expected value $1$, but instead of paying $1$ for $S_2$, I pay what I have, $S_1(y_1)$, for $S_1(y_1)S_2$.  After this second bet is settled, I have $S_1(y_1)S_2(y_2)$.%
\footnote{This argument assumes that the price of my second bet is exactly equal to my current capital after the first bet is settled.  Since the only constraint is that I not risk my capital becoming negative, we might imagine other options.  I could reserve some of my capital and buy a payoff that costs less, or perhaps I might be allowed to buy a payoff that costs more than my current capital if this payoff is bounded away from zero.  But these ideas do not really increase my betting opportunities.  When I am not allowed to risk my capital becoming negative, any bet I make can be understood as buying a nonnegative payoff that has expected value equal to my current capital.}

Multiplying betting scores may sometimes give a more reasonable evaluation of scientific exploration than other methods of combination.  Consider the scientist who uses a significance level of $5\%$ in a search for factors that might influence a phenomenon.  Her initial explorations are promising, but only after $20$ tries ($20$ slightly different chemicals in a medical study or $20$ slightly different stimuli in a psychological study) does she find an effect that is significant at $5\%$.  How seriously should we take this apparent discovery?  One standard answer is that the significance level of $5\%$ should be multiplied by $20$; this is the Bonferroni adjustment.  It has a betting rationale; we may suppose that the scientist has put up $\$1$ each time she tests a factor, thereby investing a total of $\$20$.  She loses her $\$1$ on each of the first $19$ tries, but she wins $\$20$ on her $20$th try.  When we recognize that she actually invested \$20, not merely \$1, we might conclude that her final betting score is 20/20, or 1.  But this will be unfair if the first $19$ experiments were promising, as the product of $20$ betting scores that are only a little larger than $1$ may be reasonably large.

In many fields, the increasing resources being devoted to the search for significant effects has led to widespread and justified skepticism about published statistical studies purporting to have discovered such effects.  This is true for both experimental studies and studies based on databases.  A recent replication of published experimental studies in social and cognitive psychology has shown that many of these studies cannot be replicated \cite{OSC:2015}.  A recent massive study using databases from medical practice has shown that null hypotheses known to be true are rejected at a purported $5\%$ level about $50\%$ of the time \cite{Madigan/etal:2014,Schuemie/etal:2018}.  A recent review of database studies in finance has noted that although a large number of factors affecting stock prices have been identified, few of these results seem to be believed, inasmuch as each study ignores the previous studies \cite{Harvey:2017}. These developments confirm that we need to report individual statistical results in ways that embed them into broader research programs.   Betting scores provide one tool for this undertaking, both for the individual scientist reporting on her own research and for the meta-analyst reporting on the research of a scientific community \cite{terSchure/Grunwald:2019}.

\subsection{Score for a single bet = likelihood ratio}\label{subsec:likely}

For simplicity, suppose $P$ is discrete.  Then the assumption $\mathbf{E}_P(S)=1$ can be written
\[
    \sum_y S(y)P(y)=1.  
\]
Because $S(y)$ and $P(y)$ are nonnegative for all $y$, this tells us that the product $S P$ is a probability distribution.  Write $Q$ for $S P$, and call $Q$ the alternative \emph{implied} by the bet $S$.  If we suppose further that $P(y)>0$ for all $y$, then $S=Q/P$, and 
\begin{equation}\label{eq:ratio}
      S(y) = \frac{Q(y)}{P(y)}.
\end{equation}
A betting score is a likelihood ratio.

Conversely, a likelihood ratio is a betting score.  Indeed, if $Q$ is a probability distribution for $Y$, then $Q/P$ is a bet by our definition, because $Q/P\ge 0$ and
\[
        \sum_y \frac{Q(y)}{P(y)}P(y)=\sum_y Q(y)=1. 
\]

\paragraph{When I have a hunch that $Q$ is better\ldots}

We began with your claiming that $P$ describes the phenomenon $Y$ and my making a bet $S$ satisfying $S\ge 0$ and, for simplicity, $\mathbf{E}_P(S)=1$.  There are no other constraints on my choice of $S$.  The choice may be guided by some hunch about what might work, or I may act on a whim.  I may not have any alternative distribution $Q$ in mind.  Perhaps I do not even believe that there is an alternative distribution that is valid as a description of $Y$.  

Suppose, however, that I do have an alternative $Q$ in mind.  I have a hunch that $Q$ is a valid description of $Y$.  In this case, should I use $Q/P$ as my bet?  The thought that I should is supported by Gibbs's inequality, which says that
\begin{equation}\label{eq:gibb}
   \mathbf{E}_Q\left(\ln \frac{Q}{P}\right) \ge \mathbf{E}_Q\left(\ln \frac{R}{P}\right)
\end{equation}
for any probability distribution $R$ for $Y$.%
\footnote{Many readers will recognize $\mathbf{E}_Q\left(\ln (Q/P)\right)$ as the Kullback-Leibler divergence between $Q$ and $P$.  In the terminology of Kullback's 1959 book \cite[p.~5]{Kullback:1959}, it is the mean information for discrimination in favor of $Q$ against $P$ per observation from $Q$.}  
Because any bet $S$ is of the form $R/P$ for some such $R$, \eqref{eq:gibb} tells us that $\mathbf{E}_Q(\ln S)$ is maximized over $S$ by setting $S:=Q/P$.

Why should I choose $S$ to maximize $\mathbf{E}_Q(\ln S)$?  Why not maximize $\mathbf{E}_Q(S)$?  Or perhaps $Q(S\ge 20)$ or $Q(S\ge 1/\alpha)$ for some other significance level $\alpha$?


Maximizing $\mathbf{E}(\ln S)$ makes sense in a scientific context because, as we have seen, successive betting scores are combined by multiplication.   When $S$ is the product of many successive factors, maximizing $\mathbf{E}(\ln S)$ maximizes $S$'s rate of growth.  This point was made famously and succinctly by John L.\ Kelly, Jr., in 1956 \cite{Kelly:1956}:  ``it is the logarithm which is additive in repeated bets and to which the law of large numbers applies.''  The idea has been used extensively in gambling theory \cite{Breiman:1961}, information theory \cite{Cover/Thomas:1991}, finance theory \cite{Luenberger:2014}, and machine learning \cite{Cesa-Bianchi/Lugosi:2006}.  I am proposing that we put it to greater use in statistical testing.  

We can use Kelly's insight even when betting is opportunistic and hence does not define alternative joint probabilities for successive outcomes.   Even if the null hypothesis $P$ does provide joint probabilities for a phenomenon $(Y_1,Y_2,\ldots)$, successive opportunistic bets $S_1,S_2,\ldots$ against $P$ will not imply an alternative $Q$.  Each bet $S_i$ will imply only an alternative $Q_i$ for $Y_i$ in light of the actual outcomes $y_1,\ldots,y_{n-1}$.%
\footnote{One team may test $P$ with a bet $S_1(Y_1)=Q_1(Y_1)/P(Y_1)$.  A second team, after seeing $y_1$ and the betting score $S_1(y_1)$, may then test $P$ with a bet $S_2(Y_2)=Q_2(Y_2)/P(Y_2|Y_1=y_1)$.  But because no one made a single bet $S(Y_1,Y_2)$ at the outset, no alternative probability distribution $Q$ for $(Y_1,Y_2)$ has been specified.  From $P$ and $S_1$ and $S_2$ we can obtain probabilities for $Y_1$ and conditional probabilities for $Y_2$ given $Y_1$'s observed value $y_1$, but we do not obtain conditional probabilities for $Y_2$ given other values for $Y_1$.}  
A game-theoretic law of large numbers nevertheless holds with respect to the sequence $Q_1,Q_2,\ldots$: if they are valid in the betting sense (an opponent will not multiply their capital by a large factor betting against them), then the average of the logarithms of the $Y_i$ will approximate the average of the expected values assigned by the $Q_i$.  See \cite[Chapter 2]{Shafer/Vovk:2019}.

\paragraph{The Neyman-Pearson lemma.}

In 1928 \cite{Neyman/Pearson:1928}, Jerzy Neyman and E.\ S.\ Pearson suggested that for a given significance level $\alpha$, we choose a rejection region $E$ such that $Q(y)/P(y)$ is at least as large for all $y\in E$ as for any $y\notin E$, where $Q$ is an alternative hypothesis.%
\footnote{Asking the reader's indulgence, I leave aside the difficulty that it may be impossible, especially if $P$ and $Q$ are discrete, to do this precisely or uniquely.}  
Let us call the bet $S_E$ with this choice of $E$ the \emph{level-$\alpha$ Neyman-Pearson bet} against $P$ with respect to $Q$.  The \emph{Neyman-Pearson lemma} says that this choice of $E$ maximizes 
\[
    Q(\text{test rejects $P$})=Q(Y\in E)=Q(S_E(Y)\ge 1/\alpha),
\] 
which we call the \emph{power} of the test with respect to $Q$.  In fact, $S_E$ with this choice of $E$ maximizes $Q(S(Y)\ge 1/\alpha)$ over all bets $S$, not merely over all-or-nothing bets.%
\footnote{If $S\ge 0$, $E_P(S)=1$, $0<\alpha<1$, and $Q(S\ge 1/\alpha)>0$, define an all-or-nothing bet $S'$ by
\[
  S'(y):=
  \begin{cases}
      \frac{1}{\alpha} & \text{if } S(y) \ge \frac{1}{\alpha}\\
      0                      & \text{if } S(y) < \frac{1}{\alpha}.
  \end{cases}
\]
Then $E_P(S')<1$ and $Q(S'\ge 1/\alpha)=Q(S\ge 1/\alpha)$.  Dividing $S'$ by $E_P(S')$, we obtain an all-nothing bet that has expected value $1$ under $P$ and a greater probability of exceeding $1/\alpha$ under $Q$ than $S$.} 
It does not maximize $\mathbf{E}_Q(\ln S)$ unless $Q=S_E P$, and this is usually an unreasonable choice for $Q$, because it gives probability one to $E$.

It follows from Markov's inequality that when the level-$\alpha$ Neyman-Pearson bet against $P$ with respect to $Q$ just barely succeeds, the bet $Q/P$ succeeds less:  it multiplies the money risked by a smaller factor.%
\footnote{The Neyman-Pearson bet just barely turns $1$ into $1/\alpha$ when 
               $P(Q/P \ge Q(y)/P(y))=\alpha$.
               By Markov’s inequality,
               $P(Q/P \ge Q(y)/P(y)) \le P(y)/Q(y)$.  So $1/\alpha \ge Q(y)/P(y)$.}  
But the success of the Neyman-Pearson bet may be unconvincing in such cases; see Examples 1 and~2 in Section~\ref{subsec:examples}.

R.\ A.\ Fisher famously criticized Neyman and Pearson for confusing the scientific enterprise with the problem of ``making decisions in an acceptance procedure'' \cite[Chapter 4]{Fisher:1956}.  Going beyond all-or-nothing tests to general testing by betting is a way of taking this criticism seriously.  The choice to ``reject'' or ``accept'' is imposed when we are testing a widget that is to be put on sale or returned to the factory for rework, never in either case to be tested again.  But in many cases scientists are looking for evidence against a hypothesis that may be tested again many times in many ways.

\paragraph{When the bet loses money\ldots}

In the second paragraph of the introduction, I suggested that a betting score of 5 casts enough doubt on the hypothesis being tested to merit attention.  We can elaborate on this by noting that a value of 5 or more for $S(y)$ means, according to~\eqref{eq:ratio}, that the outcome $y$ was at least 5 times as likely under the alternative hypothesis $Q$ than under the null hypothesis $P$.

Suppose we obtain an equally extreme result in the opposite direction:  $S(y)$ comes out less than $1/5$.  Does this provide enough evidence in favor of $P$ to merit attention?  Maybe and maybe not.  A low value of $S(y)$ does suggest that $P$ describes the phenomenon better than $Q$.  But $Q$ may or may not be the only plausible alternative.  It is the alternative for which the bet $S$ is optimal in a certain sense.  But as I have emphasized, we may have chosen $S$ blindly or on a whim, without any real opinion or clue as to what alternative we should consider.  In this case, the message of a low betting score is not that $P$ is supported by the evidence but that we should try a rather different bet the next time we test $P$.

This understanding of the matter accords with Fisher's contention that testing usually precedes the formulation of alternative hypotheses in science \cite[p.~246]{Bennett:1990}, \cite[p.~57]{Senn:2011}.

\begin{table}
\caption{Elements of a study that tests a probability distribution by betting.  The proposed study may be considered meritorious and perhaps even publishable regardless of its outcome when the null hypothesis $P$ and the implied alternative $Q$ are both initially plausible and the implied target is reasonably large.  A large betting score discredits the null hypothesis.}\label{ta:test}
\begin{center}
\begin{tabular}{lll}
    & name  & notation\\
\textbf{Proposed study} &&\\
\hspace*{10pt}initially unknown outcome  & phenomenon & $Y$\\
\hspace*{10pt}probability distribution for $Y$ & null hypothesis & $P$\\\addlinespace[2pt]
\hspace*{10pt}\begin{minipage}{0.4\linewidth}\raggedright
nonnegative function of $Y$ with\\ 
\hspace*{10pt}expected value $1$ under $P$\end{minipage}
                                                                     & bet & $S$\\ \addlinespace[2pt]
\hspace*{10pt}$S P$                                        & implied alternative & $Q$ \\
\hspace*{10pt}$ \exp\left(\mathbf{E}_Q(\ln S) \right)$   & implied target & $S^*$\\
\addlinespace[4pt]
\textbf{Results} &&\\
\hspace*{10pt}actual value of $Y$      & outcome & $y$ \\ \addlinespace[2pt]
\hspace*{10pt}\begin{minipage}{0.4\linewidth}\raggedright
factor by which money risked\\ 
\hspace*{10pt}has been multiplied\end{minipage}
                                           & betting score & $S(y)$
\end{tabular}
\end{center}
\end{table}

\subsection{Implied targets}\label{subsec:target}

How impressive a betting score can a scientist hope to obtain with a particular bet $S$ against $P$?  As we have seen, the choice of $S$ defines an alternative probability distribution, $Q=S P$, and $S$ is the bet against $P$ that maximizes $\mathbf{E}_Q(\ln S)$.  If the scientist who has chosen $S$ takes $Q$ seriously, then she might hope for a betting score whose logarithm is in the ballpark of $\mathbf{E}_Q(\ln S)$ --- i.e., a betting score in the ballpark of 
\[
             S^* :=  \exp\left(\mathbf{E}_Q(\ln S)\right).
\]
Let us call $S^*$ the \emph{implied target} of the bet $S$. By~\eqref{eq:gibb}, $S^*$ cannot be less than $1$.  The implied target of the all-or-nothing bet \eqref{eq:NP} is always $1/\alpha$, but as we have already noticed, that bet's implied $Q$ is not usually a reasonable hypothesis.

The notion of an implied target is analogous to Neyman and Pearson's notion of power with respect to a particular alternative.  But it has the advantage that the scientist cannot avoid discussing it by refusing to specify a particular alternative.  The implied alternative $Q$ and the implied target $S^*$ are determined as soon as the distribution $P$ and the bet $S$ are specified.  The implied target can be computed without even mentioning $Q$, because
\[
    \mathbf{E}_Q(\ln S) 
               = \sum_y Q(y)\ln S(y)
               = \sum_y P(y)S(y)\ln S(y)
               = \mathbf{E}_P(S\ln S).
\]
If bets become a standard way of testing probability distributions, the implied target will inevitably be provided by the software that implements such tests, and referees and editors will inevitably demand that it be included in any publication of results.  Even if the scientist has chosen her bet $S$ on a hunch and is not committed in any way to $Q$, it is the hypothesis under which $S$ is optimal, and a proposed test will not be interesting to others if it would not be expected to achieve much even when it was optimal.  

On the other hand, if the implied alternative is seen as reasonable and interesting in its own right, and if the implied target is high, then a proposed study may merit publication regardless of how the betting score comes out (see Table~\ref{ta:test}).  In these circumstances, even a low betting score will be informative, as it indicates that the implied alternative is no better than the null.   This feature of testing by betting may help mitigate the problem of publication bias.

\subsection{Elementary examples}\label{subsec:examples}

Many if not most misuses of p-values fall into one of three groups:
\begin{enumerate}
\item
An estimate is statistically and practically significant but hopelessly contaminated with noise.  Andrew Gelman and John Carlin contend that this case ``is central to the recent replication crisis in science'' \cite[p.~900]{Gelman/Carlin:2017}.
\item
A test with a conventional significance level and high power against a very distinct alternative rejects the null hypothesis with a borderline outcome even though the null has greater likelihood than the alternative \cite[pp.~249--250]{Dempster:1997}.  
\item
A high p-value is interpreted as evidence for the null hypothesis.  Although such an interpretation is never countenanced by theoretical statisticians, it is distressingly common in some areas of application \cite{Amrhein/etal:2019,Cready/etal:2019}.
\end{enumerate} 
To see how these three cases can be handled with betting scores, it suffices to consider elementary examples.  Here I will consider examples where the null and alternative distributions of the test statistic are normal with the same variance.


\paragraph{Example 1.}  

Suppose  $P$ says that $Y$ is normal with mean 0 and standard deviation 10, $Q$ says that $Y$ is normal with mean 1 and standard deviation 10, and we observe $y=30$.  
\begin{itemize}
\item
Statistician A simply calculates a p-value:
$
     P(Y\ge30)\approx 0.00135.
$
She concludes that $P$ is strongly discredited.
\item
Statistician B uses the Neyman-Pearson test with significance level $\alpha=0.05$, which rejects $P$ when $y>16.5$.  Its power is only about 6\%.  Seeing $y=30$, it does reject $P$.  If she used the test as a bet, the statistician has multiplied the money she risked by $20$. 
\item
Statistician C uses the bet $S$ given by
\[
       S(y):= \frac{q(y)}{p(y)}= \frac{(10\sqrt{2\pi})^{(-1)}\exp(-(y-1)^2/200)}{(10\sqrt{2\pi})^{(-1)}\exp(-y^2/200)} =\exp\left(\frac{2y-1}{200}\right), 
\]
for which
\[
     \mathbf{E}_Q(\ln(S))  =  \mathbf{E}_Q\left(\frac{2y-1}{200}\right) = \frac{1}{200},
\]
so that the implied target is $\exp(1/200)\approx 1.005$.  She does a little better than this very low target; she multiplies the money she risked by $\exp(59/200)\approx1.34$.
\end{itemize}
The power and the implied target both told us in advance that the study was a waste of time.  The betting score of $1.34$ confirms that little was accomplished, while the low p-value and the Neyman-Pearson rejection of $P$ give a misleading verdict in favor of $Q$.
\qed

\paragraph{Example 2.}  Now the case of high power and a borderline outcome:  $P$ says that $Y$ is normal with mean 0 and standard deviation 10, $Q$ says that $Y$ is normal with mean 37 and standard deviation 10, and we observe $y=16.5$.  
\begin{itemize}
\item
Statistician A again calculates a p-value:
$
    P(Y\ge16.5)\approx 0.0495.
$
She concludes that $P$ is discredited.
\item
Statistician B uses the Neyman-Pearson test that rejects when $y>16.445$.  This test has significance level $\alpha=0.05$, and its power under $Q$ is almost 98\%.  It rejects; Statistician B multiplies the money she risked by 20. 
\item
Statistician C uses the bet $S$ given by $S(y):= q(y)/p(y)$.  Calculating as in the previous example, we see that $S$'s implied target is $939$ and yet the betting score is only $S(16.5)=0.477$.  Rather than multiply her money, Statistician C has lost more than half of it.  She concludes that the evidence from her bet very mildly favors $P$ relative to $Q$.
\end{itemize}
Assuming that $Q$ is indeed a plausible alternative, the high power and high implied target suggest that the study is meritorious.  But the low p-value and the Neyman-Pearson rejection of $P$ are misleading.  The betting score points in the other direction, albeit not enough to merit attention.
\qed

\paragraph{Example 3.}  Now the case of a non-significant outcome:  $P$ says that $Y$ is normal with mean 0 and standard deviation 10, $Q$ says that $Y$ is normal with mean 20 and standard deviation 10, and we observe $y=5$. 
\begin{itemize} 
\item
Statistician A calculates the p-value
$
     P(Y\ge5)\approx 0.3085.
$
As this is not very small, she concludes that the study provides no evidence about $P$.
\item
Statistician B uses the Neyman-Pearson test that rejects when $y>16.445$.  This test has significance level $\alpha=0.05$, and its power under $Q$ is about 64\%.  It does not reject; Statistician B loses all the money she risked. 
\item
Statistician C uses the bet $S$ given by $S(y):= q(y)/p(y)$.  This time $S$'s implied target is approximately $7.39$ and yet the actual betting score is only $S(5)\approx0.368$.  Statistician C again loses more than half her money.  She again  concludes that the evidence from her bet favors $P$ relative to $Q$ but not  enough to merit attention.
\end{itemize}
In this case, the power and the implied target both suggested that the study was marginal.  The Neyman-Pearson conclusion was ``no evidence".  The bet $S$ provides the same conclusion; the score $S(y)$ favors $P$ relative to $Q$ but too weakly to merit attention.
\qed

The underlying problem in the first two examples is the mismatch between the concept of a p-value on the one hand and the concepts of a fixed significance level and power on the other.  This mismatch and the confusion it engenders disappears when we replace p-values with betting scores and power with implied target.  The bet, implied target, and betting score always tell a coherent story.  In Example~1, the implied target told us that the bet would not accomplish much, and the implied target close to 1 only confirmed this.  In Example~2, the high implied target told us that we had a good test of $P$ relative to $Q$, and $P$'s passing this test strongly suggests that $Q$ is not better than $P$.

The problem in Example~3 is the meagerness of the interpretation available for a middling to high p-value.  The theoretical statistician correctly tells us that such a p-value should be taken as ``no evidence''. But a scientist who has put great effort into a study will want to believe that its result signifies something.  In this case, the merit of the betting score is that it blocks any erroneous claim with a concrete message:  it tells us the direction the result points and how strongly.

\section{Comparing scales}\label{sec:comp}

The notion of a p-value retains a whiff of betting.  In a passage I will quote shortly, R.\ A.\ Fisher used the word ``odds'' when comparing two p-values.  But obtaining a p-value $\mathsf{p}(y)$ cannot be interpreted as multiplying money risked by $1/\mathsf{p}(y)$, because the bet that achieves this multiplication could only have been made with foreknowledge of $y$. Pretending that we had made the bet would be cheating, and some penalty for this cheating --- some sort of shrinking --- is needed to make $1/\mathsf{p}(y)$ comparable to a betting score.

The inadmissibility of $1/\mathsf{p}(y)$ as a betting score is confirmed by its infinite expected value under $P$. Shrinking it to make it comparable to a betting score means shrinking it to a payoff with expected value $1$.  In the ideal case, $\mathsf{p}(y)$ is uniformly distributed between $0$ and $1$ under $P$, and there are infinitely many ways of shrinking it to something with expected value $1$.%
\footnote{In the general case, the p-value $\mathsf{p}(y)$ given by~\eqref{eq:p2} is stochastically dominated under $P$ by the uniform distribution, and so the shrinking discussed here will produce a payoff with expected value $1$ or less.}  
Some are discussed in \cite[Section 11.5]{Shafer/Vovk:2019}.
No one has made a convincing case for any particular choice from this infinitude; the choice is fundamentally arbitrary.  But it would useful to make some such choice, because the use of p-values will never completely disappear, and if we also use betting scores, we will find ourselves wanting to compare the two scales.

It seems reasonable to shrink p-values in a way that is monontonic, smooth, and unbounded, and the exact way of doing this will sometimes be unimportant.  My favorite way, because it is easy to remember and calculate (and only for this reason), is to shrink $\mathsf{p}(y)$ to 
\begin{equation}\label{eq:adj}
                   S(y) := \frac{1}{\sqrt{\mathsf{p}(y)}}-1.
\end{equation}
Table~\ref{ta:sqp} applies this rule to some commonly used significance levels.  If we retain the conventional 5\% threshold for saying that a p-value merits attention, then this table accords with the suggestion, made in the introduction to this paper, that multiplying our money by 5 merits attention.  Multiplying our money by 2 or 3, or by 1/2 or 1/3 as in Examples~2 and~3 of Section~\ref{subsec:examples}, does not meet this threshold.

\begin{table}
\caption{Making a p-value into a betting score.}\label{ta:sqp}
\begin{center}
\begin{tabular}{C{2cm}C{2cm}C{2cm}}
p-value   & $\displaystyle\frac{1}{\text{p-value}}$   &  $\displaystyle\frac{1}{\sqrt{\text{p-value}}}-1$\\
\midrule
0.10      & 10     &  2.2\\
0.05      & 20     & 3.5\\
0.01      & 100    & 9.0\\
0.005    & 200     & 13.1\\
0.001     & 1,000    & 30.6\\
0.000001     & 1,000,000    & 999\\\bottomrule
\end{tabular}
\end{center}
\end{table}

If we adopt a standard rule for shrinking p-values, we will have a fuller picture of what we are doing when we use a conventional test that is proposed without any alternative hypothesis being specified.  Since it determines a bet, the rule for shrinking implies an alternative hypothesis.

\paragraph{Example 4.}  Consider R.\ A.\ Fisher's analysis of Weldon's dice data in the first edition of his \emph{Statistical Methods for Research Workers} \cite[pp.~66--69]{Fisher:1925}.  Weldon threw 12 dice together 26,306 times and recorded, for each throw, how many dice came up 5 or 6.  Using this data, Fisher tested the bias of the dice in two different ways.
\begin{itemize}
\item
First, he performed a $\chi^2$ goodness-of-fit test.  On none of the 26,306 throws did all 12 dice come up 5 or 6, so he pooled the outcomes 11 and 12 and performed the test with 12 categories and 11 degrees of freedom.  The $\chi^2$ statistic came out 40.748, and he noted that ``the actual chance in this case of $\chi^2$ exceeding 40$\cdot$75 if the dice had been true is $\cdot$00003.''
\item
Then he noted that in the $12\times26{,}306=315{,}672$ throws of a die there were altogether $106{,}602$ 5s and 6s.   The expected number is $315{,}672/3=105{,}224$ with standard error 264.9, so that the observed number exceeded expectation by 5.20 times its standard error, and ``a normal deviation only exceeds $5{\cdot2}$ times its standard error once in 5 million times.''  
\end{itemize}
Why is the one p-value so much less than the other?  Fisher explained:  ``The reason why this last test gives so much higher odds than the test for goodness of fit, is that the latter is testing for discrepancies of any kind, such, for example, as copying errors would introduce.  The actual discrepancy is almost wholly due to a single item, namely, the value of $p$, and when that point is tested separately its significance is more clearly brought out.'' (Here $p$ is the probability of a 5 or 6, hypothesized to be 1/3.)

The transformation~\eqref{eq:adj} turns the p-values 0.00003 and 1 in 5 million into betting scores (to one significant figure) 200 and 2,000, respectively.  This does not add much by itself, but it brings a question to the surface.  The statistician has chosen particular tests and could have chosen differently. What alternative hypotheses are implied when the tests chosen are considered as bets?  

For simplicity, consider Fisher's second test and the normal approximation he used.  With this approximation, the frequency
\[
                 Y := \frac{\text{total number of 5s and 6s}}{315{,}672}
\]
is normally distributed under the null hypothesis $P$, with mean $1/3$ and standard deviation $0.00084$.  The observed value $y$ is $106{,}602/315{,}672\approx0.3377$.  As Fisher noted, the deviation from $1/3$, $0.0044$, is $5.2$ times the standard deviation.   The function $\mathsf{p}(y)$ for Fisher's test is 
\[
       \mathsf{p}(y) = 2\left( 1 - \Phi\left( \frac{\lvert y-\frac13\rvert}{0.00084} \right) \right),
\]
where $\Phi$ is the cumulative distribution function for the standard normal distribution.  The density $q$ for the alternative $Q$, obtained by multiplying $P$'s normal density $p$ by~\eqref{eq:adj} is symmetric around $1/3$, just as $p$ is.  It has the same value at $1/3$ as $p$ does, but much heavier tails.  These heavy tails approximately double the standard deviation.  The probability of a deviation of $0.0044$ or more under $Q$ is still very small, but only about 1 in a thousand instead of 1 in 5 million.

A different rule for shrinking the p-value to a betting score will of course produce a different alternative hypothesis $Q$.  But a wide range of rules will give roughly the same picture.

We can obtain an alternative hypothesis in the same way for the $\chi^2$ test.  Whereas the distribution of the $\chi^2$ statistic is approximately normal under the null hypothesis, the alternative will have much heavier multi-dimensional tails.  The existence of such an alternative, even if vaguely defined, supports Joseph Berkson's classic argument for discretion when using such a test \cite{Berkson:1938}.

\section{Betting games as statistical models}\label{sec:estimation}

In the preceding section, we learned how a single probability distribution can be tested by betting.  In this section, we look at how this mode of testing extends to testing composite hypotheses and estimating parameters. 

The extension will be obvious to anyone familiar with how standard tests are extended from point to composite hypotheses and used to form confidence sets.  A composite hypothesis is rejected if each of its elements is rejected, and a $(1-\alpha)$-confidence set consists of all hypotheses not rejected at level $\alpha$.  But when we test by betting, it is easy to get confused about who is doing the betting, and so a greater degree of formality is helpful.  This formality can be provided by the notion of a \emph{testing protocol}, which is studied in great detail and used as a foundation for mathematical probability in \cite{Shafer/Vovk:2019}.  A testing protocol tells how betting offers are made, who decides what offers to accept, and who decides the outcomes.  It is then the protocol, not a probability distribution or a parametric class of probability distributions, that represents the phenomenon.  

According to R.\ A.\ Fisher \cite{Fisher:1922}, the theory of statistical estimation begins with the assumption that the statistician has only partial knowledge of a probability distribution describing a phenomenon. %
She knows only that this probability distribution is in a known class $(P_\theta)_{\theta\in\Theta}$. The corresponding assumption in the betting picture is that the statistician stands outside a testing protocol, instructing the bettor to play a certain strategy but seeing only some of the moves.  The parameter $\theta$ is one move she does not see.  But if she believes that the protocol is a valid description of the phenomenon and has no reason to think that a betting strategy for Skeptic that she has specified has been exceptionally lucky, she can rely on the presumption that it will not multiply the capital it risks by a large factor to claim \emph{warranties} that resemble the direct probability statements made by 19th-century statisticians \cite{Shafer:2019} and confidence intervals as defined by Jerzy Neyman in the 1930s \cite{Neyman:1937}.

This section shows how testing protocols can represent statistical models and used to derive warranties.  Section~\ref{subsec:testing} introduces protocols for testing a single probability distribution.  Section~\ref{subsec:para} introduces protocols for testing parametric statistical models and defines precisely the notion of a $(1/\alpha)$-warranty.  Section~\ref{subsec:least} discusses how these ideas apply to non-parametric estimation by least squares.

\subsection{Testing protocols}\label{subsec:testing}

The following protocol formalizes Section~\ref{sec:test}'s method of testing a probability distribution $P$ for a phenomenon $Y$ that takes values in a set $\mathcal{Y}$. 

\medskip

\indentII Skeptic announces $S:\mathcal{Y}\to[0,\infty)$ such that $\mathbf{E}_P(S)= 1$.\\
\indentII Reality announces $y\in\mathcal{Y}$.\\
\indentII $\K := S(y)$.

\medskip

\noindent
Like all testing protocols considered in this paper, this is a perfect information protocol; the players move sequentially and each sees the other's move as it is made.

This representation generalizes naturally to successive variables $Y_1,\ldots,Y_N$ each purported to be described by $P$ and purported to be mutually independent:

\medskip

\indentI $\K_0:=1$.\\
\indentI FOR $n=1,2,\dots,N$:\\
\indentII Skeptic announces $S_n:\mathcal{Y}\to[0,\infty)$ such that $\mathbf{E}_P(S_n)= \K_{n-1}$.\\
\indentII Reality announces $y_n\in\mathcal{Y}$.\\
\indentII $\K_n := S_n(y_n)$.

\medskip

\noindent
For any payoff $S:\mathcal{Y}^N\to[0,\infty)$ such that $\mathbf{E}_{P^N}(S)=1$, Skeptic can play so that $\K_N=S(y_1,\ldots,y_N)$; on the $n$th round, he makes the bet $S_n$ given by $S_n(y):=\mathbf{E}_{P^N}(S|y_1,\ldots,y_{n-1},y)$.

We can similarly use an $N$-round protocol to test any hypothesized discrete stochastic process $Y_1,\ldots,Y_N$ --- i.e., any probability distribution $P$ on $\mathcal{Y}^N$.  In this case, the condition $\mathbf{E}_P(S_n)= \K_{n-1}$ on Skeptic's $n$th-round bet is replaced by $\mathbf{E}_P(S_n\vert y_1,\ldots,y_{n-1})= \K_{n-1}$.

Little or nothing is gained merely by representing a bet on an $N$-dimensional outcome as $N$ successive bets, but the $N$-round representation allows further generalizations.  One very important generalization is the introduction of a signal $x_n$ announced by Reality at the beginning of each round.  When such signals are given, preceding signals can be used along with preceding outcomes to price Skeptic's moves on each round.  As part of the protocol, we specify a probability distribution $P_{x_1,y_1\ldots,x_{n-1},y_{n-1},x_{n}}$ for each round $n$ and each possible sequence of signals and outcomes $x_1,y_1\ldots,x_{n-1},y_{n-1},x_{n}$ that might precede Skeptic's move on that round; the condition on Skeptic's move $S_n$ is then that $\mathbf{E}_{P_{x_1,y_1\ldots,x_{n-1},y_{n-1},x_{n}}}=\K_{n-1}$.  The now-classical way of handling signals or ``independent variables'' $x_1,\ldots,x_N$, going back at least to Fisher \cite{Aldrich:2005b}, is to assume or pretend that they are constants, fixed and known at the outset.  Here that pretense is not needed. 

Another level of generality is attained when we introduce a third player, Forecaster, who announces the probability distribution for $y_n$ just before Skeptic makes his $n$th round bet.  Then we are testing Forecaster, not necessarily a particular probability distribution or system of probability distributions.  In this general picture, a distribution or system of distributions is a strategy for Forecaster in the protocol.

Protocols with signals and Forecasters are studied in \cite{Shafer/Vovk:2019}.  I will not discuss them further in this paper, but in Section~\ref{subsec:least} we will see an example of another generalization, where Skeptic is offered a smaller menu of bets on each round.  This generalization takes us into the domain that statisticians call non-parametric and that others have dubbed imprecise probability \cite{Augustin/etal:2014}.

\subsection{Parametric models as testing protocols}\label{subsec:para}

Suppose $(P_\theta)_{\theta\in\Theta}$ is a parametric model for a variable $Y$ that takes values in the set $\mathcal{Y}$, and consider the following testing protocol: 

\medskip

\indentII Reality announces $\theta\in\Theta$.\\
\indentII Skeptic announces $S:\mathcal{Y}\to[0,\infty)$ such that $E_{P_\theta}(S)= 1$.\\
\indentII Reality announces $y\in\mathcal{Y}$.\\
\indentII $\K := S(y)$.

\medskip

\noindent
Here Skeptic is testing Reality's assertion that $P_\theta$ describes $Y$; he begins with unit capital, and $\K$ is his final capital.  

A strategy $\mathcal{S}$ for Skeptic in this protocol assigns a bet, say $\mathcal{S}^\theta$, to every $\theta\in\Theta$.  We suppose that the statistician does not see Reality's move $\theta$ but has chosen a strategy $\mathcal{S}$ for Skeptic and directed him to play it.  She does observe $y$ and then calculates Skeptic's final capital $\K$ as a function of $\theta$, say $\K_\mathcal{S}(\theta):=\mathcal{S}^\theta(y)$.  This allows her to test hypotheses and state warranties:
\begin{itemize}
\item  
For each $\theta\in\Theta$, $\K_\mathcal{S}(\theta)$ is \emph{$\mathcal{S}$'s betting score against the hypothesis $\theta$}.  
\item
For each subset $\Theta_0 \subseteq \Theta$,  
$
   \inf \{\K_\mathcal{S}(\theta)\st\theta\in\Theta_0\}
$
is \emph{$\mathcal{S}$'s betting score against the composite hypothesis $\Theta_0$}.  
\item
For each $\alpha>0$, \emph{$\mathcal{S}$ provides a $(1/\alpha)$-warranty} for the subset $W_{1/\alpha}^\mathcal{S}$ of $\Theta$ given by
\[
           W_{1/\alpha}^\mathcal{S}:= \{\theta\in\Theta \st \K_\mathcal{S}(\theta) < \frac{1}{\alpha}\}.
\]
Skeptic multiplied the capital he risked by at least $1/\alpha$ if $\theta$ is not in $W_{1/\alpha}^\mathcal{S}$.
\end{itemize}
This notion of a $(1/\alpha)$-warranty has roots in the work of Claus-Peter Schnorr \cite{Schnorr:1971} and Leonid Levin \cite{Levin:1976}.

For small $\alpha$, the statistician will tend to believe that the true $\theta$ is in $W_{1/\alpha}^\mathcal{S}$.  For example, she will not expect Skeptic to have multiplied his capital by 1000 and hence will believe that $\theta$ is in $W_{1000}^\mathcal{S}$.  But this belief is not irrefutable.  If she sees other evidence that $\theta$ is not in $W_{1000}^\mathcal{S}$, then she may conclude that Skeptic actually did multiply his capital by $1000$ using $\mathcal{S}$.  See \cite{Fraser/etal:2018} for examples of outcomes that cast doubt on confidence statements and would also cast doubt on warranties.

Every $(1-\alpha)$-confidence set has a $(1/\alpha)$-warranty.  This is because a $(1-\alpha)$-confidence set is specified by testing each $\theta$ at level $\alpha$; the $(1-\alpha)$-confidence set consists of the $\theta$ not rejected.  When $\mathcal{S}$ makes the all-or-nothing bet against $\theta$ corresponding to the test used to form the confidence set, $K_\mathcal{S}(\theta)<1/\alpha$ if and only if $\theta$ was not rejected, and hence $W_{1/\alpha}^\mathcal{S}$ is equal to the confidence set.    

Warranty sets are nested:  $W_{1/\alpha}^\mathcal{S}\subseteq W_{{1/\alpha'}^\mathcal{S}}$ when $\alpha>\alpha'$.  Standard statistical theory also allows nesting; sets with different levels of confidence can be nested.  But the different confidence sets will be based on different tests \cite{Cox:1958,Xie/Singh:2013}.  The $(1/\alpha)$-warranty sets for different $\alpha$ all come from the same strategy for Skeptic.

The notion of a $(1/\alpha)$-warranty extends to testing protocols with multiple rounds.  It can be used, for example, in the following protocol, which represents a parametric model $(P_\theta)_{\theta\in\Theta}$ for independent and identically distributed random variables $Y_1,Y_2,\ldots,Y_N$ drawn from a set $\mathcal{Y}$.

\medskip
 
\indentI $\K_0:=1$.\\
\indentI Reality announces $\theta\in\Theta$.\\
  \indentI FOR $n=1,2,\dots,N$:\\
    \indentII Skeptic announces $S_n:\mathcal{Y}\to[0,\infty)$ such that $\mathbf{E}_{P_\theta}(S_n) =\K_{n-1}$.\\
    \indentII Reality announces $y_n\in\mathcal{Y}$.\\
    \indentII $\K_n := S_n(y_n)$.

\medskip

\noindent
Here $\K_0$ is Skeptic's initial capital, $S_n$ is his bet on the $n$th round, and $\K_n$ is his capital at the end of the $n$th round.  In this protocol a strategy $\mathcal{S}$ for Skeptic is a more complicated object, but it still makes the final capital $\K_N$ a function of $\theta$ and hence determines $(1/\alpha)$-warranties. 

If $T:\Theta\times\mathcal{Y}^N\to[0,\infty)$ satisfies $\mathbf{E}_{P_\theta}(T(\theta,\cdot)) =1$ for all $\theta\in\Theta$, then Skeptic has a strategy $\mathcal{S}$ with final capital $\K_\mathcal{S}(\theta)$ equal to $T(\theta,y_1,\ldots,y_N)$.  (See the explanation immediately following the protocol with $N$ rounds in Section~\ref{subsec:testing}.)  In particular, for every event $E$ in the protocol (every measurable subset $E$ of $\mathcal{Y}^N$), Skeptic has a strategy $\mathcal{S}$ that produces the capital
\[   
    \K_\mathcal{S}(\theta)
    = 
       \begin{cases}
           \frac{1}{P_\theta^N(E)} & \text{if $y_1,\ldots,y_n\in E$}\\
           0 & \text{otherwise}.     
      \end{cases}  
\]
It follows that any $(1-\alpha)$-confidence set for $\theta$ is a $1/\alpha$-warranty set for $\theta$.

Instead of stopping the protocol after some fixed number $N$ of rounds, we might suppose that it continues indefinitely and think about the sequence of $(1/\alpha)$-warranty sets obtained from a continuing strategy.  The intuition supporting the $(1/\alpha)$-warranty set obtained from an initial sequence of rounds will still apply to the $(1/\alpha)$-warranty set obtained from a longer sequence; Skeptic simply continued betting.  This contrasts with the theory of confidence sets based on sequences of observations; the tests used for a confidence interval based on $N$ trials are usually not used for the confidence interval based on a longer sequence, say $N'>N$ of trials, and so using both intervals raises issues about multiple testing and ``sampling to a foregone conclusion'' \cite{Cornfield:1966,Shafer/etal:2011}.

How should the statistician choose the strategy for Skeptic?  An obvious goal is to obtain small warranty sets.  But a strategy that produces the smallest warranty set for one $N$ and one warranty level $1/\alpha$ will not generally do so for other values of these parameters.  So any choice will be a balancing act.  How to perform this balancing act is an important topic for further research; see \cite{Grunwald/etal:2019}.

\subsection{Non-parametric least squares}\label{subsec:least}

Consider the following protocol, which represents in a minimal way assumptions about errors of measurement suggested by Gauss in 1821 \cite{Stigler:1986}:  each error is bounded, and an error of a given size might just as well be negative as positive.

\medskip

  \indentI $\K_0 :=1$.\\
  \indentI FOR $n=1,2,\ldots$:\\
    \indentII Skeptic announces $M_n \in [-\K_{n-1},\K_{n-1}]$.\\
    \indentII Reality announces $e_n\in[-1,1]$.\\
    \indentII $\K_n := \K_{n-1} + M_n e_n$.
 
\medskip

\noindent
On the $n$th round, Skeptic can buy a multiple of the error $e_n$ at the price $0$.   The protocol does not specify the number $N$ of rounds; play will stop at some point, but we need not specify a stopping time in advance.

Were $P$ a probability distribution for $e_n$ that gives probability one to the interval $[-1,1]$ and has mean zero, then we could write $\mathbf{E}_P(\K_{n-1} + M_n e_n)=\K_{n-1}$.  If we allowed Skeptic to choose any nonnegative payoff that has expected value $\K_{n-1}$ under $P$, even if it is not of the form $\K_{n-1} + M_n e_n$, and if we also specified a stopping time, then the protocol would be a special case of general type described at the end of Section~\ref{subsec:testing}. But we have specified neither a probability distribution $P$ nor a stopping time $N$.

Although no probability distribution is specified, the protocol can be taken as theory and tested just as a protocol representing a probability distribution is tested.  As a theory about actual errors of $e_1,e_2,\ldots$, the protocol asserts that Skeptic will not multiply his capital by a large factor.  It is discredited to the extent that Skeptic succeeds in doing so.

In \cite[Section 3.3]{Shafer/Vovk:2019}, it is shown that Skeptic has a strategy, based on Hoeff\-ding's inequality, that guarantees $\K_n \ge 20$ for every $n$ such that $\vert \overline{e}_n \vert > 2.72\sqrt{n}$.  If Skeptic plays this strategy and eventually reaches an $n$ for which $\vert \overline{e}_n \vert > 2.72\sqrt{n}$, then he can claim a betting score of $20$ against the protocol.  Had we specified a probability distribution $P$ on $(0,1)$ with mean zero and allowed Skeptic to choose on each round any nonnegative payoff with expected value under $P$ equal to his current capital, then he could have chosen a particular value $N$ large enough for the central limit theorem to be effective and claimed a betting score of $20$ if $\vert \overline{e}_N \vert > 2.0\sqrt{N}$.  But this would work only for particular value $N$.

Now suppose that the errors are added to a quantity $\mu$ that is being measured.  As in Section~\ref{subsec:para}, assume that the statistician stands outside the game and sees only the measurements $y_1,y_2\ldots$.  She does not see the quantity $\mu$ being measured or the errors $e_1,e_2,\ldots$. Then we have this protocol.

\medskip

\indentI $\K_0:=1$.\\
\indentI Reality announces $\mu\in\bbbr$.\\
  \indentI FOR $n=1,2,\ldots$:\\
    \indentII Skeptic announces $M_n\in[-\K_{n-1},\K_{n-1}]$.\\
    \indentII Reality announces $e_n\in[-1,1]$ 
      and sets $y_n:=\mu+ e_n$.\\
    \indentII $\K_n := \K_{n-1} + M_n e_n$.

\medskip

Now the strategy for Skeptic that guarantees $\K_n \ge 20$ for every $n$ such that $\vert \overline{e}_n \vert > 2.72\sqrt{n}$ can be used to obtain warranties for $\mu$.  After $100$ measurements, for example, it gives a 20-warranty that $\mu$ is in $\overline{y}_{100} \pm 0.272$.  

The statistician will know the betting score that Skeptic has achieved only as a function of $\mu$.  But a meta-analyst, imagining that Skeptic has used his winnings from each study in the next study, can multiply the functions of $\mu$ to obtain warranties about $\mu$ that may be more informative than those from the individual studies.

Averaging measurements of a single quantity to estimate the quantity measured is the most elementary instance of estimation by least squares.  The ideas developed here extend to the general theory of estimation by least squares,  in which $\mu$ is multi-dimensional and multi-dimensional signals $x_1,x_2,\ldots$ are used.  An asymptotic theory with this generality is developed in \cite[Section 10.4]{Shafer/Vovk:2019}.

\section{Conclusion}\label{sec:conclusion}

The probability calculus began as a theory about betting, and its logic remains the logic of betting, even when it serves to describe phenomena.  But in their quest for the appearance of objectivity, mathematicians have created a language (likelihood, significance, power, p-value, confidence) that pushes betting into the background.  

This deceptively objective statistical language can encourage overconfidence in the results of statistical testing and neglect of relevant information about how the results are obtained.  In recent decades this problem has become increasingly salient, especially in medicine and the social sciences, as numerous influential statistical studies in these fields have turned out to be misleading.

In 2016, the American Statistical Association issued a statement listing common misunderstandings of p-values and urging full reporting of searches that produce p-values \cite{Wasserstein/Lazar:2016}.  Many statisticians fear, however, that the situation will not improve.  Most dispiriting are studies showing that both teachers of statistics and scientists who use statistics are apt to answer questions about the meaning of p-values incorrectly \cite{McShane/Gal:2017,Gigerenzer:2018}.  Andrew Gelman and John Carlin argue persuasively that the most frequently proposed solutions (better exposition, confidence intervals instead of tests, practical instead of statistical significance, Bayesian interpretation of one-sided p-values, and Bayes factors) will not work \cite{Gelman/Carlin:2017}.  The only solution, they contend, is ``to move toward a greater acceptance of uncertainty and embracing of variation'' (p.~901).

In this context, the language of betting emerges as an important tool of communication.  When statistical tests and conclusions are framed as bets, everyone understands their limitations.  Great success in betting against probabilities is the best evidence we can have that the probabilities are wrong, but everyone understands that such success may be mere luck.  Moreover, candor about the betting aspects of scientific exploration can communicate truths about the games scientists must and do play --- honest games that are essential to the advancement of knowledge.

This paper has developed new ways of expressing statistical results with betting language.  The basic concepts are \emph{bet} (not necessarily all-or-nothing), \emph{betting score} (equivalent to likelihood ratio when the bets offered define a probability distribution), \emph{implied target} (an alternative to power), and \emph{$(1/\alpha)$-warranty} (an alternative to $(1-\alpha)$-confidence).  Substantial research is needed to apply these concepts to complex models, but their greatest utility may be in communicating the uncertainty of simple tests and estimates.

\section{Acknowledgements}

Many conversations over the past several years have inspired and influenced this paper.  Most important, perhaps, were conversations about game-theoretic testing and meta-analysis with Peter Gr\"unwald and Judith ter Schure at the Centrum Wiskunde \& Informatica in Amsterdam in December 2018.  Also especially important were conversations with Gert de Cooman and Jasper Bock and their students at the University of Ghent, with Harry Crane, Jacob Feldman, Robin Gong, Barry Loewer, and others in Rutgers University's seminar on the Foundations of Probability, and with Jason Klusowski, William Strawderman, and Min-ge Xie in the seminar of the Statistics Department at Rutgers.

Many others have provided useful feedback after this paper was first drafted, including John Aldrich,  Jasper De Bock,  Steve Goodman, Prakash Gorroochurn, Sander Greenland, Alan H\'ajek,  Wouter Koolen, David Madigan, Deborah Mayo, Teddy Seidenfeld, Stephen Senn, Nozer Singpurwalla, Mike Smithson, Aris Spanos, Matthias Troffaes, Conor Mayo-Wilson, Vladimir Vovk, and several anonymous referees.  My understanding of the idea of testing by trying to multiply the money one risks has developed over several decades in the course of my collaboration with Vovk.

\section{Appendix: Situating testing by betting}\label{sec:sit}

\paragraph{Disclaimers.}  My theme has been that we can communicate statistical conclusions and their uncertainty more effectively with betting scores than with p-values.  This is not to say that well trained statisticians are unable to use p-values effectively.  We have been using them effectively for two centuries. 

I have emphasized that testing by betting is not a chapter in decision theory, because tests can use amounts of money so small that no one cares about them.  The betting is merely to make a point, and play money would do.  This is not to say that decision theory is an unimportant chapter in statistical methodology.  Many statistical problems do involve decisions for which the utilities and probabilities required by various decision theories are available.  These theories include the Neyman-Pearson theory and Bayesian theory.  These theories do not use p-values, and so replacing p-values by betting scores would not affect them.

\paragraph{Are the probabilities tested subjective or objective?}  The probabilities may represent someone's opinion, but the hypothesis that they say something true about the world is inherent in the project of testing them.

\paragraph{Are the probabilities tested frequencies?}  Not in any concrete sense.  In the general protocols of Section~\ref{subsec:testing}, Skeptic can select any particular event to which $P$ assigns a probability and adopt a strategy that will produce a large betting score unless Reality makes the frequency of the event approximate that probability.  But only the most salient probabilities will be tested, and as Abraham Wald pointed out in the 1930s, only a countable number of them could be tested \cite{Bienvenu/etal:2009}. So the identification of $P$'s probabilities with frequencies is always approximate and usually hypothetical.  The connection with frequencies is even more tenuous when the theory tested involves limited betting offers, as in non-parametric and other imprecise-probability models.

\paragraph{How is testing by betting related to Bruno de Finetti's theory of probability?} In de Finetti's picture, an individual who has \emph{evidence about facts} expresses the strength of that evidence by \emph{odds he offers} and prices at which he would buy or sell \cite{de-Finetti:1970}.  Testing by betting uses \emph{how a bet came out} as a measure of \emph{evidence about probabilities}.  These are different undertakings.

The contrast is especially clear in protocols in which a forecaster chooses and announces betting offers.  De Finetti emphasized the viewpoint of the forecaster. Testing by betting emphasizes the viewpoint of Skeptic, who decides how to bet, and the perspective of the statistician who considers the evidence provided by the outcomes.  Most authors on imprecise probabilities, including Peter Walley \cite{Walley:1991}, have also emphasized the forecaster's viewpoint.

Some readers have asked how the expectation that Skeptic not obtain a large betting score is related to de Finetti's condition of coherence.  The two are not closely related.  Coherence is a condition on offers made to Skeptic:  they should not allow him to make money no matter what Reality does.  This condition is met by all the protocols in this paper.  Usually Reality can move so that Skeptic obtains a large betting score, and hence the expectation that this not happen becomes, when a strategy for Skeptic is fixed, an expectation about Reality’s moves.  If Reality violates this expectation, we doubt whether the protocol is a valid description of Reality.  When an unknown parameter is involved, the belief that the protocol is valid tells us something about the parameter.

\paragraph{Why is the proposal to test by betting better than other proposals for remedying the misuse of p-values?}  Many authors have proposed remedying the misuse of p-values by supplementing them with additional information.%
\footnote{Some idea of the range of proposals is provided by \cite{Wasserstein/etal:2019} and \cite{Mayo:2018}.}    
Sometimes the additional information involves Bayesian calculations \cite{Bayarri/etal:2016,Matthews:2018}. Sometimes it involves likelihood ratios \cite{Colquhoun:2019}.  Sometimes it involves attained power \cite{Mayo/Spanos:2006}.

I find nearly all these proposals persuasive as ways of correcting the misunderstandings and misinterpretations to which p-values are susceptible.  Each of them might be used by a highly trained mathematical statistician to explain what has gone wrong to another highly trained mathematical statistician.  But adding more complexity to the already overly complex idea of a p-value may not help those who are not specialists in mathematical statistics.  We need strategies for communicating with millions of people.  Worldwide, we teach p-values to millions every year, and hundreds of thousands of them may eventually use statistical tests in one way or another.  

The most important argument for betting scores as a replacement for p-values is its simplicity.  I do not know any other proposal that is equally simple.

\paragraph{Is testing by betting a novel idea?}

The idea of testing purported probabilities by betting is part of our culture.  But its first appearance in mathematics seems to be in the work of Jean Ville in 1939 \cite[pp.~87--89]{Ville:1939}. Richard von Mises's principle of the impossibility of a gambling system \cite[p.\ 25]{vonMises:1928} said that a strategy for selecting throws on which to bet should not improve a gambler's chances.  Ville replaced this with a stronger and more precise condition:  the winnings from a strategy that risks only a fixed amount of money should not increase indefinitely, even if the strategy is allowed to vary the amount bet and the side on which to bet.  As Ville put it, a player who follows such a strategy ``ne gagne pas ind\'efiniment'' \cite[p.~89]{Ville:1939}.

Nevertheless, testing by betting remains remarkably absent in mathematical statistics. Some authors allude to betting from time to time, but they usually mention only all-or-nothing bets.  Betting never takes center stage, even though it is central to everyone's intuitions about probability.  Ever since Jacob Bernoulli, we have downplayed the betting aspect of testing and estimation in order to make them look more objective --- more scientific.

Many mathematicians, including R.\ A.\ Fisher, have advocated using likelihood as a direct measure of evidence. 
In his \emph{Statistical Methods and Scientific Inference} \cite{Fisher:1956}, Fisher suggested that in most cases the plausible values of a parameter indexing a class of probability distributions are those for which the likelihood is at least $(1/15)$th of its maximum value.%
\footnote{On pp.~75--75 of the third edition, which appeared in 1973, Fisher used diagrams to show ``outside what limits the likelihood falls to levels at which the corresponding values of the parameter become implausible.''  In these diagrams, ``zones are indicated showing the limits within which the likelihood exceeds 1/2, 1/5, and 1/15 of the maximum.  Values of the parameter outside the last limit are obviously open to grave suspicion.''}
Later authors, including A.\ W.\ F.\ Edwards \cite{Edwards:1972} and Richard Royall \cite{Royall:1997}, published book-length arguments for using likelihood ratios to measure the strength of evidence, and articles supporting this viewpoint continue to appear.  
I have not found in this literature any allusion to the idea that a likelihood ratio measures the success of a bet against a null hypothesis.

Because Ville introduced martingales into probability theory to test probabilities, and because the notion was further developed in this direction by Per Martin-L\"of and Claus-Peter Schnorr \cite{Bienvenu/etal:2009}, we might expect that statisticians who use martingales would recognize the possibility of interpreting nonnegative martingales directly as tests.  They know that when $P$ and $Q$ are probability distributions for a sequence $Y_1,Y_2,\ldots$, the sequence of likelihood ratios
\begin{equation}\label{eq:mart}
    1,\frac{Q(Y_1)}{P(Y_1)},\frac{Q(Y_1,Y_2)}{P(Y_1,Y_2)},\ldots
\end{equation}
is a nonnegative martingale, and they would see no novelty in the following observations:
\begin{itemize}
\item 
If I am allowed to bet on $Y_1,\ldots,Y_n$ at rates given by $P$, then I can buy the payoff $Q(Y_1,\ldots,Y_n)/P(Y_1,\ldots,Y_n)$ for one monetary unit. 
\item
If I bet sequentially, on each $Y_k$ at rates given by $P(Y_k\vert y_1,\ldots,y_{k-1})$ when  $y_1,\ldots,y_{k-1}$ are known, and I always use my winnings so far to buy that many units of $Q(Y_k\vert y_1,\ldots,y_{k-1})/P(Y_k\vert y_1,\ldots,y_{k-1})$, then~\eqref{eq:mart} will be the capital process resulting from my strategy.
\end{itemize}
But they also know that Joseph L.\ Doob purified the notion of a martingale of its betting heritage when he made it a technical term in modern probability theory.  Martingales are now widely used in sequential analysis, time series, and survival analysis \cite{Aalen/etal:2009,Lai:2009}, but I have not found in the statistical literature any use of the idea that the successive realized values of a martingale already represent, without being translated into the language of p-values and significance levels, the cumulative evidential value of the outcomes of a betting strategy.

\addcontentsline{toc}{section}{References}

\bibliographystyle{plain}
   \InputIfFileExists{54.bbl}


\end{document}